\documentclass{article}
\usepackage{amssymb}
\usepackage{amsfonts}
\usepackage{amsmath}

\begin{document}

\begin{center}
{\Large About special elements in quaternion algebras\\
 over finite fields}

\begin{equation*}
\end{equation*}%
Diana SAVIN 
\begin{equation*}
\end{equation*}
\end{center}

\textbf{Abstract. }{\small In this paper we study special Fibonacci quaternions and special generalized Fibonacci-Lucas quaternions in quaternion algebras over finite fields.}


\bigskip

\begin{equation*}
\end{equation*}

\textbf{1. Introduction}%
\begin{equation*}
\end{equation*}

Let $K$ be an algebraic number field and let $\mathbb{H}_{K}\left( \alpha ,\beta \right) ~$\textit{be the
generalized quaternion algebra} over the field $K.$ An arbitrary element from this algebra has the form $x=x_{1}\cdot 1+x_{2}e_{2}+x_{3}e_{3}+x_{4}e_{4},$
where $x_{i}\in K,i\in \{1,2,3,4\}$, where the elements of the basis $%
\{1,e_{2},e_{3},e_{4}\}$ satisfy the following multiplication table: \vspace{%
3mm}

\begin{center}
\begin{tabular}{ccccc}
$\cdot \,\,$\vline & $1$ & $e_{2}$ & $e_{3}$ & $e_{4}$ \\ \hline
$1$\thinspace \vline & $1$ & $e_{2}$ & $e_{3}$ & $e_{4}$ \\ 
$e_{2}$\vline & $e_{2}$ & $\alpha $ & $e_{4}$ & $\alpha e_{3}$ \\ 
$e_{3}$\vline & $e_{3}$ & $-e_{4}$ & $\beta $ & $-\beta e_{2}$ \\ 
$e_{4}$\vline & $e_{4}$ & $-\alpha e_{3}$ & $\beta e_{2}$ & $-\alpha \beta $%
\end{tabular}
\end{center}

\bigskip

Let $x$$\in$$\mathbb{H}_{K}\left( \alpha ,\beta \right) $ and let $\boldsymbol{n}\left(x\right) $ be the norm of a generalized
quaternion $x.$ This norm has the following expression $\boldsymbol{n}\left(
x\right) =x_{1}^{2}-\alpha x_{2}^{2}-\beta x_{3}^{2}+\alpha \beta x_{4}^{2}.$%
We recall that  $\mathbb{H}_{K}\left( \alpha ,\beta \right)$ is a division algebra if and only if for $x\in \mathbb{H}%
_{K}\left( \alpha ,\beta \right) $ we have $\boldsymbol{n}\left( x\right) =0$
if and only if $x=0.\,$ Otherwise, the algebra $\mathbb{H}_{K}\left( \alpha
,\beta \right) $ is called \textit{a split }algebra. \\
When $\alpha=\beta=-1$ and $K=\mathbb{R}$ we obtain Hamilton quaternion algebra $\mathbb{H}_{\mathbb{R}}\left( -1 ,-1 \right),$ with the basis $%
\{1,i,j,k\}.$ It is known that this algebra is a division algebra.\\
When $\alpha=\beta=-1$ and $K$ is a finite field ($K=\mathbb{Z}_{p},$ where $p$ is a prime positive integer) we obtain the quaternion algebra $\mathbb{H}_{\mathbb{Z}_{p}}\left( -1 ,-1 \right).$ It is known that this algebra is a split algebra ([Gr, Mi, Ma; 15], Theorem 2, p. 882).\\

\medskip

Let $\left( f_{n}\right) _{n\geq 0}$ be the Fibonacci sequence: 
\begin{equation*}
f_{0}=0;f_{1}=1;f_{n}=f_{n-1}+f_{n-2},\;n\geq 2
\end{equation*}%
and $\left( l_{n}\right) _{n\geq 0}$ be the Lucas sequence: 
\begin{equation*}
l_{0}=2;l_{1}=1;l_{n}=l_{n-1}+l_{n-2},\;n\geq 2
\end{equation*}%
In his paper [Ho; 63], A. F. Horadam introduced the Fibonacci quaternions $\left( F_{n}\right) _{n\geq 0}:$ 
\begin{equation*}
F_{n}=f_{n}+f_{n+1}\cdot i+f_{n+2}\cdot j+f_{n+3}\cdot k, \;n\geq 0.
\end{equation*} 
\begin{equation*}
\end{equation*}
Some properties of Fibonacci quaternions in certain quaternion algebras were investigated in
the papers [Ak, Ko, To; 14], [Ta; Ya; 15], [Po, Ki, Ke; 16], [Y$\ddot{u}$, Ay; 16], etc.
In this paper we determine the Fibonacci quaternions (respectively generalized Fibonacci-
Lucas quaternions) which are zero divisors in the quaternion ring $\mathbb{H}_{\mathbb{Z}_{p}}\left( -1 ,-1 \right)$ and the Fibonacci quaternions which are invertible in the quaternion ring $\mathbb{H}_{\mathbb{Z}_{p}}\left( -1 ,-1 \right).$\\
We recall some necessary results about Fibonacci sequence and about Fibonacci quaternions.\\
\smallskip\\
\textbf{Proposition 1.1} ([Fib]). \textit{Let} $\left( f_{n}\right) _{n\geq 0}$ \textit{be the Fibonacci sequence. If} $f_{n}$ \textit{is a prime number} ($n\geq5$), \textit{then} $n$
\textit{is a prime number.}\\
\smallskip\\
\textbf{Proposition 1.2} ([Fib]). \textit{Let} $\left( f_{n}\right) _{n\geq 0}$ \textit{be the Fibonacci sequence. If} $n,m$$\in$$\mathbb{N}^{*},$ $n|m,$ \textit{then} $f_{n}|f_{m}.$\\
\smallskip\\
\textbf{Theorem 1.1} ([Su], [Su; 92]).
\textit{Let} $p\notin \left\{2, 5\right\}$ \textit{be a prime positive integer. Then}
\[
f_{\frac{p-\left(\frac{p}{5}\right)}{2}}\equiv   
\left\{\begin{array}{ll}
           \ 0\: (mod\:  p), if\: p\equiv1 (mod 4),\\
           \ 2\left(-1\right)^{\left[\frac{p+5}{10}\right]}\cdot \left(\frac{p}{5}\right)\cdot 5^{\frac{p-3}{4}}\: (mod\:  p) , $ if $ \: p\equiv 3 (mod 4)
           \end {array}
\right.
\]
\textit{where} $\left[x\right]$ \textit{is the integer part of x and} $\left(\frac{p}{5}\right)$ \textit{is the Legendre symbol.}\\
\smallskip\\
\textbf{Proposition 1.3} ([Fib]). \textit{The cycle of the Fibonacci numbers mod} $2$ \textit{is}
$$0,1,1 (0,1,1),...$$
so:\\
i) \textit{the cycle-length of the Fibonacci numbers mod} $2$ \textit{is} $3;$\\
ii) \textit{a Fibonacci number} $f_{n}$ \textit{is even if and only if} $n$$\equiv$$0$ (\textit{mod} $3$).\\
iii) \textit{a Lucas number} $l_{n}$ \textit{is even if and only if} $n$$\equiv$$0$ (\textit{mod} $3$).\\
\smallskip\\
\textbf{Proposition 1.4} ([Fib]). \textit{The cycle of the Fibonacci numbers mod} $3$ \textit{is}
$$0,1,1,2,0,2,2,1 (0,1,1,2,0,2,2,1),... $$
\textit{so, the cycle-length of the Fibonacci numbers mod} $3$ \textit{is} $8.$\\
\smallskip\\
\textbf{Proposition 1.5} ([Fib]). \textit{The cycle of the Fibonacci numbers mod} $5$ \textit{is}
$$0,1,1,2,3,0,3,3,1,4,0,4,4,3,2,0,2,2,4,1 (0,1,1,2,3,0,3,3,1,4,0,4,4,3,2,0,2,2,4,1),...$$
so: \\
i) \textit{the cycle-length of the Fibonacci numbers mod} $5$ \textit{is} $20;$\\
ii) $f_{n}$$\equiv$$0$ (\textit{mod} $5$) \textit{if and only if} $n$$\equiv$$0$ (\textit{mod} $5$).\\
\smallskip\\
\textbf{Proposition 1.6} ([Ho; 63]). \textit{Let} $\left( f_{n}\right) _{n\geq 0}$ \textit{be the Fibonacci sequence and} $F_{n} (n\geq 0)$ \textit{be the Fibonacci quaternions. Then, the norm of the nth Fibonacci 
quaternions is:}
$$n\left(F_{n}\right)=3f_{2n+3}.$$
\smallskip\\
Let $p$ be an odd prime positive integer. In the paper [Mi, Se; 11], C. J. Miguel, R. Serodio calculated the number of zero-divisors in $\mathbb{H}_{\mathbb{Z}_{p}}\left( -1 ,-1 \right).$
They proved that:
\bigskip\\
\textbf{Theorem 1.2} ([Mi, Se; 11]). \textit{The number of zero-divisors in} $\mathbb{H}_{\mathbb{Z}_{p}}\left( -1 ,-1 \right)$ \textit{is}:
$$p^{3}+p^{2}-p.$$
\smallskip\\
Let $p$ be a prime positive integer. According to [Da, Dr; 70], if $f_{z}$ is the smallest
non-zero Fibonacci number with the property $p|f_{z},$ then $z=z\left(p\right)$ is defined as the entry point of $p$ in the Fibonacci sequence.\\
\smallskip\\
\textbf{Remark 1.1.} ([Da, Dr; 70]). $p|f_{n}$ \textit{if and only if} $z\left(p\right)|n,$ $\left(\forall\right) n\in \mathbb{N}.$\\
\smallskip\\
\textbf{Definition 1.1.} ([Da, Dr; 70]). \textit{Let} $\left( f_{n}\right) _{n\geq 0}$ \textit{be the Fibonacci sequence. If we reduce the Fibonacci sequence modulo} $p,$ \textit{to obtain a periodic sequence. The period} $k=k\left(p\right)$ \textit{is the smallest integer} $k$ \textit{for which}
$$f_{k}\equiv0\ \left(mod\ p\right)\ and\ f_{k+1}\equiv1\ \left(mod\ p\right).$$
\smallskip\\
For example: $k\left(2\right)=z\left(2\right)=3,$ $k\left(3\right)=8,$ $z\left(3\right)=4,$ $k\left(5\right)=20,$ $z\left(5\right)=5$ (according to Proposition 1.3, Proposition 1.4, Proposition 1.5).\\
\smallskip\\
\textbf{Remark 1.2.} ([Da, Dr; 70]). \textit{The entry point} $z\left(p\right)$ \textit{divides the period} $k\left(p\right).$\\
\smallskip\\
\textbf{Remark 1.3.} ([Da, Dr; 70]). \textit{Let} $p$ \textit{be an odd prime positive integer. The
following statements are true:}\\
i) $k\left(p\right)=z\left(p\right)$ \textit{if} $z\left(p\right)\equiv 2$ (\textit{mod} $4$);\\
ii) $k\left(p\right)=2z\left(p\right)$ \textit{if} $z\left(p\right)\equiv 0$ (\textit{mod} $4$);\\
iii) i) $k\left(p\right)=4z\left(p\right)$ \textit{if} $z\left(p\right)$ \textit{is odd}.\\
\smallskip\\
\textbf{Proposition 1.7.} ([Fib]).\textit{Let} $(f_{n})_{n\geq 0}$ \textit{be the
Fibonacci sequence} and $(l_{n})_{n\geq 0}$ \textit{be the
Lucas sequence} \textit{Then:}\newline
i)  \begin{equation*}
f^{2}_{n}+f^{2}_{n+1}=f_{2n+1},\forall ~n\in \mathbb{N};
\end{equation*}%
ii)
\begin{equation*}
f_{n}+f_{n+4}=3f_{n+2},\forall ~n\in \mathbb{N};
\end{equation*}%
iii)
\begin{equation*}
l^{2}_{n}+l^{2}_{n+1}=5f_{2n+1},\forall ~n\in \mathbb{N};
\end{equation*}%
iv)
\begin{equation*}
f_{m}l_{m+p}=f_{2m+p} + \left(-1\right)^{m+1}\cdot f_{p},\forall ~m,p\in \mathbb{N};
\end{equation*}%
v)
\begin{equation*}
f_{n}+f_{n+2}=l_{n+1},\forall ~n\in \mathbb{N};
\end{equation*}%

\begin{equation*}
\end{equation*}
\textbf{2. Special Fibonacci quaternions in the quaternion algebra $\mathbb{H}_{\mathbb{Z}_{p}}\left( -1 ,-1 \right).$}%
\begin{equation*}
\end{equation*}
In this section we determine the zero divisors Fibonacci quaternions, respectively the invertible Fibonacci quaternions in the quaternion rings $\mathbb{H}_{\mathbb{Z}_{2}}\left( -1 ,-1 \right),$ $\mathbb{H}_{\mathbb{Z}_{3}}\left( -1 ,-1 \right)$, $\mathbb{H}_{\mathbb{Z}_{5}}\left( -1 ,-1 \right),$ $\mathbb{H}_{\mathbb{Z}_{p}}\left( -1 ,-1 \right),$ where $p$ is a prime positive integer, $p>5.$ The idea for to find these elements is the following: since $\mathbb{H}_{\mathbb{Z}_{p}}\left( -1 ,-1 \right),$ (where $p$ is a prime positive integer) is a split algebra, it is isomorphic to the matrix ring $\textit{M}_{2}\left(\mathbb{Z}_{p}\right).$ If 
$\varphi$ $:$ $\mathbb{H}_{\mathbb{Z}_{p}}\left( -1 ,-1 \right)$$\rightarrow$$\textit{M}_{2}\left(\mathbb{Z}_{p}\right)$ is a rings isomorphism, then a Fibonacci quaternion $F_{n}$ is a zero divisor in $\mathbb{H}_{\mathbb{Z}_{p}}\left( -1 ,-1 \right)$ if and only if $\varphi\left(F_{n}\right)$ is a zero divisor in $\textit{M}_{2}\left(\mathbb{Z}_{p}\right)$ 
(see [Mi, Se; 11], Theorem 2.1). Moreover, it is know that a matrix $A$$\in$$\textit{M}_{2}\left(\mathbb{Z}_{p}\right)$ is a zero divisor in $\textit{M}_{2}\left(\mathbb{Z}_{p}\right)$ if and only if it is not invertible in $\textit{M}_{2}\left(\mathbb{Z}_{p}\right).$ These involve that, a Fibonacci quaternion is a zero divisor in $\mathbb{H}_{\mathbb{Z}_{p}}\left( -1 ,-1 \right)$ if and only if it is not invertible in $\mathbb{H}_{\mathbb{Z}_{p}}\left( -1 ,-1 \right).$ So, a Fibonacci quaternion $F_{n}$ is a zero divisor in $\mathbb{H}_{\mathbb{Z}_{p}}\left( -1 ,-1 \right)$ if and only if $n\left(F_{n}\right)=0.$\\
\bigskip\\
\textbf{Proposition 2.1.} i) \textit{Fibonacci quaternions} $F_{n}$ \textit{are zero divisors in the quaternion ring} $\mathbb{H}_{\mathbb{Z}_{2}}\left( -1 ,-1 \right)$ \textit{if and only if} $n$$\equiv$$0$ (\textit{mod} $3$).\\
ii) \textit{In} $\mathbb{H}_{\mathbb{Z}_{2}}\left( -1 ,-1 \right)$ \textit{there are one Fibonacci quaternion which is zero divisor and two Fibonacci quaternions which are invertible.}
\smallskip\\
\textbf{Proof.} i) According to Proposition 1.6 and Proposition 1.3 we obtain that Fibonacci quaternions $F_{n}$ are zero divisors in the quaternion ring $\mathbb{H}_{\mathbb{Z}_{2}}\left( -1 ,-1 \right)$ if and only if $f_{2n+3}$$\equiv$$0$ (mod $2$) if and only if $n$$\equiv$$0$ (mod $3$).\\
ii) There are $16$ elements in the quaternion algebra  $\mathbb{H}_{\mathbb{Z}_{2}}\left( -1 ,-1 \right)$. According to Theorem 1.2 we know that there are $10$ zero divisors in $\mathbb{H}_{\mathbb{Z}_{2}}\left( -1 ,-1 \right).$
 We want to determine how many of them are zero divisors Fibonacci quaternions and after this we determine how many invertible Fibonacci quaternions are in the quaternion algebra 
$\mathbb{H}_{\mathbb{Z}_{2}}\left( -1 ,-1 \right).$\\
Applying i) and Proposition 1.3 we obtain that there is one Fibonacci quaternion which is zero divisor in $\mathbb{H}_{\mathbb{Z}_{2}}\left( -1 ,-1 \right),$ namely for $f_{n}=\overline{0},$ $f_{n+1}=\overline{1},$ $f_{n+2}=\overline{1},$ $f_{n+3}=\overline{0}$ (in 
$\mathbb{Z}_{2}$), so $F_{n}=i+j.$\\ 
Applying i) and Proposition 1.3, it results immediately that there are two invertible Fibonacci quaternions in $\mathbb{H}_{\mathbb{Z}_{2}}\left( -1 ,-1 \right)$ : $F_{n_{1}}=\overline{1}+i+k$ and $F_{n_{2}}=\overline{1}+j+k.$
\smallskip\\
\textbf{Proposition 2.2.} \textit{All Fibonacci quaternions} $F_{n}$ \textit{are zero divisors in the quaternion ring} $\mathbb{H}_{\mathbb{Z}_{3}}\left( -1 ,-1 \right)$ \textit{and the number of these elements is} $8.$\\
\smallskip\\
\textbf{Proof.} Applying Proposition 1.6 we obtain that $n\left(F_{n}\right)=\overline{0}$ in $\mathbb{Z}_{3},$ $\left(\forall\right) n$$\in$$\mathbb{N}.$ This implies that all Fibonacci quaternions $F_{n}$ are zero divisors in the quaternion ring $\mathbb{H}_{\mathbb{Z}_{3}}\left( -1 ,-1 \right).$\\
There are $81$ elements in the quaternion algebra  $\mathbb{H}_{\mathbb{Z}_{2}}\left( -1 ,-1 \right).$ According to Theorem 1.2, $33$ of them are zero divisors. Now, we determine how many of these $33$ elements are zero divisors Fibonacci quaternions. From the previously proved
and Proposition 1.4 we obtain the following zero divisors Fibonacci quaternions in $\mathbb{H}_{\mathbb{Z}_{3}}\left( -1 ,-1 \right):$\\
$F_{n_{1}}=i+j+\overline{2}k,$ $F_{n_{2}}=\overline{1}+i+\overline{2}j,$ $F_{n_{3}}=\overline{1}+\overline{2}i+\overline{2}k,$ $F_{n_{4}}=\overline{2}+\overline{2}j+\overline{2}k,$        $F_{n_{5}}=\overline{2}i+\overline{2}j+k,$
$F_{n_{6}}=\overline{2}+\overline{2}i+j,$   $F_{n_{7}}=\overline{2}+i+k,$ $F_{n_{8}}=\overline{1}+j+k.$\\
\smallskip\\
\textbf{Proposition 2.3.} i) \textit{Fibonacci quaternions} $F_{n}$ \textit{are zero divisors in the quaternion ring} $\mathbb{H}_{\mathbb{Z}_{5}}\left( -1 ,-1 \right)$ \textit{if and only if} $n$$\equiv$$1$ (\textit{mod} $5$).\\
ii) \textit{In} $\mathbb{H}_{\mathbb{Z}_{5}}\left( -1 ,-1 \right)$ \textit{there are} $4$ \textit{Fibonacci quaternions which are zero divisors and} $16$ \textit{Fibonacci quaternions which are invertible.}
\smallskip\\
\textbf{Proof.} i) According to Proposition 1.6 and Proposition 1.5 (ii) we obtain that Fibonacci quaternions $F_{n}$ are zero divisors in the quaternion ring $\mathbb{H}_{\mathbb{Z}_{5}}\left( -1 ,-1 \right)$ if and only if $f_{2n+3}$$\equiv$$0$ ( mod $5$) if and only if $2n+3$$\equiv$$0$ (mod $5$) if and only if $n$$\equiv$$1$ (mod $5$).\\
ii) There are $625$ elements in the quaternion algebra  $\mathbb{H}_{\mathbb{Z}_{2}}\left( -1 ,-1 \right).$ According to Theorem 1.2, $145$ of them are zero divisors. Applying i) and Proposition 1.5  we obtain $4$ zero divisor Fibonacci quaternion in $\mathbb{H}_{\mathbb{Z}_{5}}\left( -1 ,-1 \right):$\\
$F_{n_{1}}=\overline{1}+i+\overline{2}j+\overline{3}k,$ $F_{n_{2}}=\overline{3}+\overline{3}i+j+\overline{4}k,$ $F_{n_{3}}=\overline{4}+\overline{4}i+\overline{3}j+\overline{2}k,$ 
$F_{n_{4}}=\overline{2}+\overline{2}i+\overline{4}j+k.$\\
Applying i) and Proposition 1.5, it results immediately that there are $16$ invertible Fibonacci quaternions in $\mathbb{H}_{\mathbb{Z}_{5}}\left( -1 ,-1 \right)$ :
$F_{n_{1}}=i+j+\overline{2}k,$ $F_{n_{2}}=\overline{1}+\overline{2}i+\overline{3}j,$
$F_{n_{3}}=\overline{2}+\overline{3}i+\overline{3}k,$ $F_{n_{4}}=\overline{3}+\overline{3}j+\overline{3}k,$ $F_{n_{5}}=\overline{3}i+\overline{3}j+k,$ $F_{n_{6}}=\overline{3}+i+\overline{4}j,$ $F_{n_{7}}=\overline{1}+\overline{4}i+\overline{4}k,$ $F_{n_{8}}=\overline{4}+\overline{4}j+\overline{4}k,$ $F_{n_{9}}=\overline{4}i+\overline{4}j+\overline{3}k,$ $F_{n_{10}}=\overline{4}+\overline{3}i+\overline{2}j,$ $F_{n_{11}}=\overline{3}+\overline{2}i+\overline{2}k,$ $F_{n_{12}}=\overline{2}+\overline{2}j+\overline{2}k,$ $F_{n_{13}}=\overline{2}i+\overline{2}j+\overline{4}k,$ $F_{n_{14}}=\overline{2}+\overline{4}i+j,$ $F_{n_{15}}=\overline{4}+i+k,$
 $F_{n_{16}}=\overline{1}+j+k.$\\
\smallskip\\
We asked our self what happens in the general case, when we work in the quaternion ring $\mathbb{H}_{\mathbb{Z}_{p}}\left( -1 ,-1 \right).$ We obtain the following result:
\smallskip\\
\textbf{Theorem 2.1.} \textit{Let} $p$ \textit{be a prime positive integer,} $p$$\equiv$$13,17$ (mod $20$), $l$ \textit{be a positive integer and} $n=\frac{p\cdot\left(2l+1\right)+2l-1}{4}-1.$ \textit{Then, the Fibonacci quaternions} $F_{n}$ \textit{are zero divisors in the quaternion ring} $\mathbb{H}_{\mathbb{Z}_{p}}\left( -1 ,-1 \right).$\\
\smallskip\\
\textbf{Proof.} We remark that $4$$|$$\left(p\cdot\left(2l+1\right)+2l-1\right),$ so $n$ is 
 an integer number.\\
Since $p$$\equiv$$13,17$ (mod $20$), it results that the Legendre symbol $\left(\frac{p}{5}\right)=-1$ and $\frac{p-\left(\frac{p}{5}\right)}{2}=\frac{p+1}{2}$ is an odd integer number. 
Since $n=\frac{p\cdot\left(2l+1\right)+2l-1}{4}-1,$ it results $\frac{p+1}{2}|\left(2n+3\right).$ Applying Proposition 1.2, it results that
$f_{\frac{p+1}{2}}|f_{2n+3}.$ According to Theorem 1.1 and Proposition 1.6, we obtain that $n\left(F_{n}\right)=\overline{0}$ in $\mathbb{Z}_{p}.$ So, if $p$ is a prime positive integer, $p$$\equiv$$13,17$ (mod $20$) and $n=\frac{p\cdot\left(2l+1\right)+2l-1}{4}-1,$ $l$$\in$$\mathbb{N}^{*},$ then the Fibonacci quaternions $F_{n}$ are zero divisors in the quaternion ring $\mathbb{H}_{\mathbb{Z}_{p}}\left( -1 ,-1 \right).$\\
\smallskip\\
Now, we generalize the results from Proposition 2.1 and Proposition 2.3, obtaining a necessary and sufficient condition for that a Fibonacci quaternion $F_{n}$ to be a zero divisor in the quaternion ring $\mathbb{H}_{\mathbb{Z}_{p}}\left( -1 ,-1 \right).$ The difference between  the expression of $n$ from Theorem 2.1 and the expression of $n$ from Theorem 2.2 is that $n$ from Theorem 2.2 depends of $z\left(p\right).$
\smallskip\\
\textbf{Theorem 2.2.} \textit{Let} $p$ \textit{be a prime positive integer,} $p>5.$ \textit{Then}:\\
i) \textit{Fibonacci quaternions} $F_{n}$ \textit{are zero divisors in the quaternion ring} $\mathbb{H}_{\mathbb{Z}_{p}}\left( -1 ,-1 \right)$ \textit{if and only if} $n=\frac{\left(2l+1\right)\cdot z\left(p\right)-1}{2}-1,$ \textit{where} $l$$\in$$\mathbb{N}.$\\
ii) \textit{In the quaternion ring} $\mathbb{H}_{\mathbb{Z}_{p}}\left( -1 ,-1 \right)$ \textit{there are} $4$ \textit{Fibonacci quaternions which are zero divisors and} $k\left(p\right)-4$ \textit{Fibonacci quaternions which are invertible.}\\
\smallskip\\
\textbf{Proof.} i) Let $F_{n}$ be a Fibonacci quaternion which is zero divisor in the quaternion ring $\mathbb{H}_{\mathbb{Z}_{p}}\left( -1 ,-1 \right).$ Applying Proposition 1.6, it results that $p|f_{2n+3}.$ This is equivalent with $z\left(p\right)|\left(2n+3\right)$ (according to Remark 1.1). This is equivalent with $z\left(p\right)$ is an odd number and $\left(\exists\right)$ $l$$\in$$\mathbb{N}$ such that $2n+3=\left(2l+1\right)\cdot z\left(p\right)$. We obtain that 
$n=\frac{\left(2l+1\right)\cdot z\left(p\right)-1}{2}-1,$ $l$$\in$$\mathbb{N}.$\\
ii) Since $z\left(p\right)$ is an odd number, using Remark 1.3 and a reasoning similar to the one used in Proposition 2.1 (ii) and Proposition 2.3 (ii), we obtain that Fibonacci quaternions $F_{n}$ are zero divisors in the quaternion ring $\mathbb{H}_{\mathbb{Z}_{p}}\left( -1 ,-1 \right)$ if and only if $f_{2n+3}$$\equiv$$0$ ( mod $p$) if and only if $2n+3$$\equiv$$0$ (mod $z\left(p\right)$)(according to Remark 1.1) if and only if $2n$$\equiv$$z\left(p\right)-3$ (mod $z\left(p\right)$). Since $g.c.d.\left(2, z\left(p\right)\right)=1,$ it results that $\overline{2}$ is invertible in $\mathbb{Z}_{z\left(p\right)},$ so, the congruence $2n$$\equiv$$z\left(p\right)-3$ (mod $z\left(p\right)$) has unique solution mod $z\left(p\right).$ Applying Remark 1.3, we get that in $\mathbb{H}_{\mathbb{Z}_{p}}\left( -1 ,-1 \right)$ there are $4$ Fibonacci quaternions which are zero divisors and $k\left(p\right)-4$ Fibonacci quaternions which are invertible.\\
\begin{equation*}
\end{equation*}

\textbf{3. Generalized Fibonacci-Lucas quaternions which are zero divisors in the quaternion ring $\mathbb{H}_{\mathbb{Z}_{r}}\left( -1 ,-1 \right).$}%
\begin{equation*}
\end{equation*}

Let $\left( f_{n}\right) _{n\geq 0}$ be the Fibonacci sequence and $\left( l_{n}\right) _{n\geq 0}$ be the Lucas sequence.

In the paper [Fl, Sa; 15], we introduced the generalized Fibonacci-Lucas numbers and the generalized  Fibonacci-Lucas quaternions $\left( F_{n}\right) _{n\geq 0},$ namely:\\
let $n,p,q$ be arbitrary integers, $n\geq0$. The numbers $\left(g^{p,q}_{n}\right) _{n\geq 0},$
$$g^{p,q}_{n+1}=pf_{n}+ql_{n+1}, \ n\geq0$$
are called \textit{the generalized Fibonacci-Lucas numbers}.

Let $r$ be a prime positive integer and let $\mathbb{H}_{\mathbb{Z}_{r}}\left( -1 ,-1 \right)$
be the quaternion algebra over the finite field $\mathbb{Z}_{r}$ with a basis $\left\{1,i,j,k\right\},$ where $i^{2}=-1,$ $j^{2}=-1,$ $k=ij=-ji.$ The element

\begin{equation*}
G^{p,q}_{n}=g^{p,q}_{n}+g^{p,q}_{n+1}\cdot i+g^{p,q}_{n+2}\cdot j+g^{p,q}_{n+3}\cdot k
\end{equation*} 
is called \textit{the} $n$-\textit{th generalized Fibonacci-Lucas quaternion} (in the quaternion algebra $\mathbb{H}_{\mathbb{Z}_{r}}\left( -1 ,-1 \right)$).\\
\smallskip\\
We determine the generalized  Fibonacci-Lucas quaternions which are zero divisors in the quaternion ring $\mathbb{H}_{\mathbb{Z}_{r}}\left( -1 ,-1 \right),$ when $r=2$ or $r=3$ or $r=5$ or $r=p$ or $r=q.$\\
First, we compute the norm for the $n$-th generalized Fibonacci-Lucas quaternion in $\mathbb{H}_{\mathbb{Z}_{r}}\left( -1 ,-1 \right).$\\
\smallskip\\
\textbf{Proposition 3.1.} \textit{Let} $n,p,q$ \textit{be arbitrary integers}, $n\geq0$.
\textit{Let} $G^{p,q}_{n}$ \textit{be the} $n$-\textit{th generalized Fibonacci-Lucas quaternion. Let} $r$ \textit{be a prime positive
integer. Then, the norm of} $G^{p,q}_{n}$ \textit{in the quaternion algebra} $\mathbb{H}_{\mathbb{Z}_{r}}\left( -1 ,-1 \right),$ \textit{is given by}
$$\textbf{n}\left(G^{p,q}_{n}\right)=\overline{3p^{2}f_{2n+1}+15q^{2}f_{2n+3}+6pq\left(f_{2n+1}+f_{2n+3}\right)}\ (in\ \mathbb{Z}_{r}).$$
\textbf{Proof.} $$n\left(G^{p,q}_{n}\right)=G^{p,q}_{n}\overline{G^{p,q}_{n}}=\left(g^{p,q}_{n}\right)^{2}+\left(g^{p,q}_{n+1}\right)^{2}+\left(g^{p,q}_{n+2}\right)^{2}+\left(g^{p,q}_{n+3}\right)^{2}=$$
$$=\left(pf_{n-1}+ql_{n}\right)^{2} + \left(pf_{n}+ql_{n+1}\right)^{2}+\left(pf_{n+1}+ql_{n+2}\right)^{2}+\left(pf_{n+2}+ql_{n+3}\right)^{2}=$$
$$=p^{2}\left(f^{2}_{n-1}+f^{2}_{n}+f^{2}_{n+1}+f^{2}_{n+2}\right)+
q^{2}\left(l^{2}_{n}+l^{2}_{n+1}+l^{2}_{n+2}+l^{2}_{n+3}\right)+$$
$$+2pq\left(f_{n-1}l_{n}+f_{n}l_{n+1}+f_{n+1}l_{n+2}+f_{n+2}l_{n+3}\right).$$
Applying Proposition 1.7 (i, iii, iv), we have:
$$n\left(G^{p,q}_{n}\right)=p^{2}\left(f_{2n-1}+f_{2n+3}\right)+q^{2}\left(5f_{2n+1}+5f_{2n+5}\right)+$$
$$+ 2pq\left[f_{2n-1}+\left(-1\right)^{n}\cdot f_{1}+ f_{2n+1}+\left(-1\right)^{n+1}\cdot f_{1} + f_{2n+3}+\left(-1\right)^{n+2}\cdot f_{1} + f_{2n+5}+\left(-1\right)^{n+3}\cdot f_{1}\right].$$
Using Proposition 1.7 (ii), we obtain:
$$\textbf{n}\left(G^{p,q}_{n}\right)=\overline{3p^{2}f_{2n+1}+15q^{2}f_{2n+3}+6pq\left(f_{2n+1}+f_{2n+3}\right)}\ (in\ \mathbb{Z}_{r}).$$
\smallskip\\
\textbf{Proposition 3.2.} \textit{Let} $p,q$ \textit{be two 
integers} \textit{and let} $n$ \textit{be an arbitrary positive integer. Then, the following statements are true:}\\
i) \textit{if} $p$$\equiv$$0$ (\textit{mod} $2$) \textit{and} $q$$\equiv$$1$ (\textit{mod} $2$), \textit{then the generalized Fibonacci-Lucas quaternions} $G^{p,q}_{n}$ \textit{are zero divisors in the quaternion ring} $\mathbb{H}_{\mathbb{Z}_{2}}\left( -1 ,-1 \right)$ \textit{if and only if} $n$$\equiv$$0$ (\textit{mod} $3$);\\
ii) \textit{if} $p$$\equiv$$1$ (\textit{mod} $2$) \textit{and} $q$$\equiv$$0$ (\textit{mod} $2$) \textit{then the generalized Fibonacci-Lucas quaternions} $G^{p,q}_{n}$ \textit{are zero divisors in the quaternion ring} $\mathbb{H}_{\mathbb{Z}_{2}}\left( -1 ,-1 \right)$ \textit{if and only if} $n$$\equiv$$1$ (\textit{mod} $3$);\\
iii) \textit{if} $p; q$$\equiv$$0$ (\textit{mod} $2$) \textit{then all generalized Fibonacci-Lucas quaternions} $G^{p,q}_{n}$ \textit{are zero divisors in the quaternion ring} $\mathbb{H}_{\mathbb{Z}_{2}}\left( -1 ,-1 \right);$\\
iv) \textit{if} $p; q$$\equiv$$1$ (\textit{mod} $2$), \textit{then the generalized Fibonacci-Lucas quaternions} $G^{p,q}_{n}$ \textit{are zero divisors in the quaternion ring} $\mathbb{H}_{\mathbb{Z}_{2}}\left( -1 ,-1 \right)$ \textit{if and only if} $n$$\equiv$$2$ (\textit{mod} $3$);\\
\smallskip\\
\textbf{Proof.} i) The generalized Fibonacci-Lucas quaternions $G^{p,q}_{n}$ are zero divisors in the quaternion ring $\mathbb{H}_{\mathbb{Z}_{2}}\left( -1 ,-1 \right)$ if and only if $n\left(G^{p,q}_{n}\right)=\overline{0}$ (in $\mathbb{Z}_{2}$ ). Applying Proposition 3.1 and the facts that $r=2$ and $p$$\equiv$$0$ (mod $2$) and $q$$\equiv$$1$ (mod $2$), it results that
$n\left(G^{p,q}_{n}\right)=\overline{f_{2n+3}}$ in $\mathbb{Z}_{2}.$ Applying Proposition 1.3, $\overline{f_{2n+3}}=\overline{0}$ in $\mathbb{Z}_{2}$ if and only if $n$$\equiv$$0$ (mod $3$).\\
ii) If $p$$\equiv$$1$ (mod $2$) and $q$$\equiv$$0$ (mod $2$), similarly as in i) we get that the generalized Fibonacci-Lucas quaternions $G^{p,q}_{n}$ are zero divisors in the quaternion ring $\mathbb{H}_{\mathbb{Z}_{2}}\left( -1 ,-1 \right)$ if and only if $\overline{f_{2n+1}}=\overline{0}$ in $\mathbb{Z}_{2}$ if and only if $n$$\equiv$$1$ (mod $3$).\\
iii) If $p; q$$\equiv$$0$ (mod $2$), it results that $n\left(G^{p,q}_{n}\right)=\overline{0}$ (in $\mathbb{Z}_{2}$) for $\left(\forall\right)$ $G^{p,q}_{n},$ so, we obtain immediate that all generalized Fibonacci-Lucas quaternions $G^{p,q}_{n}$ are zero divisors in the quaternion ring $\mathbb{H}_{\mathbb{Z}_{2}}\left( -1 ,-1 \right).$\\
iv) If $p; q$$\equiv$$1$ (mod $2$), applying Proposition 1.7 (v) it results that $n\left(G^{p,q}_{n}\right)=\overline{f_{2n+1}}+\overline{f_{2n+3}}=\overline{l_{2n+2}}$ (in $\mathbb{Z}_{2}$ ). Using a reasoning as in i) and Proposition 1.3 (iii), we obtain that the generalized Fibonacci-Lucas quaternions $G^{p,q}_{n}$ are zero divisors in the quaternion ring $\mathbb{H}_{\mathbb{Z}_{2}}\left( -1 ,-1 \right)$ if and only if $\overline{l_{2n+2}}=\overline{0}$ (in $\mathbb{Z}_{2}$) if and only if $n$$\equiv$$2$ (mod $3$).\\
\smallskip\\
\textbf{Proposition 3.3.} \textit{Let} $p,q$ \textit{be two 
integers} \textit{and let} $n$ \textit{be an arbitrary positive integer. Then, all generalized Fibonacci-Lucas quaternions} $G^{p,q}_{n}$ \textit{are zero divisors in the quaternion ring} $\mathbb{H}_{\mathbb{Z}_{3}}\left( -1 ,-1 \right).$\\
\smallskip\\
\textbf{Proof.} Applying Proposition 3.1 we obtain that $n\left(G^{p,q}_{n}\right)=\overline{0}$ in $\mathbb{Z}_{3},$ $\left(\forall\right) n$$\in$$\mathbb{N}.$ This implies that all generalized Fibonacci-Lucas quaternions $G^{p,q}_{n}$ are zero divisors in the quaternion ring $\mathbb{H}_{\mathbb{Z}_{3}}\left( -1 ,-1 \right).$\\
\smallskip\\
\textbf{Proposition 3.4.} \textit{Let} $p,q$ \textit{be two 
integers} \textit{and let} $n$ \textit{be an arbitrary positive integer. Then, the following statements are true:}\\
i) \textit{if} $p$$\equiv$$0$ (\textit{mod} $5$), \textit{then all generalized Fibonacci-Lucas quaternions} $G^{p,q}_{n}$ \textit{are zero divisors in the quaternion ring} $\mathbb{H}_{\mathbb{Z}_{5}}\left( -1 ,-1 \right)$;\\
ii) \textit{if} $5$ \textit{does not divide} $p$ and $p$$\equiv$$q$ (\textit{mod} $5$), \textit{then the generalized Fibonacci-Lucas quaternions} $G^{p,q}_{n}$ \textit{are zero divisors in the quaternion ring} $\mathbb{H}_{\mathbb{Z}_{5}}\left( -1 ,-1 \right)$ \textit{if and only if} $n$$\equiv$$4$ (\textit{mod} $5$);\\
iii) \textit{if} $5$ \textit{does not divide} $p$ and $p+q$$\equiv$$0$ (\textit{mod} $5$), \textit{then the generalized Fibonacci-Lucas quaternions} $G^{p,q}_{n}$ \textit{are zero divisors in the quaternion ring} $\mathbb{H}_{\mathbb{Z}_{5}}\left( -1 ,-1 \right)$ \textit{if and only if} $n$$\equiv$$0$ (\textit{mod} $5$);\\
iv) \textit{if} $p$$\equiv$$1$ (\textit{mod} $5$), $q$$\equiv$$2$ (\textit{mod} $5$) \textit{or} \textit{if} $p$$\equiv$$3$ (\textit{mod} $5$), $q$$\equiv$$1$ (\textit{mod} $5$) \textit{or} \textit{if} $p$$\equiv$$2$ (\textit{mod} $5$), $q$$\equiv$$4$ (\textit{mod} $5$) \textit{or} \textit{if} $p$$\equiv$$4$ (\textit{mod} $5$), $q$$\equiv$$3$ (\textit{mod} $5$),\textit{then the generalized Fibonacci-Lucas quaternions} $G^{p,q}_{n}$ \textit{are zero divisors in the quaternion ring} $\mathbb{H}_{\mathbb{Z}_{5}}\left( -1 ,-1 \right)$ \textit{if and only if} $n$$\equiv$$1$ (\textit{mod} $5$);\\
v) \textit{if} $p$$\equiv$$2$ (\textit{mod} $5$), $q$$\equiv$$1$ (\textit{mod} $5$) \textit{or} \textit{if} $p$$\equiv$$3$ (\textit{mod} $5$), $q$$\equiv$$4$ (\textit{mod} $5$)
\textit{or} \textit{if}  $p$$\equiv$$4$ (\textit{mod} $5$), $q$$\equiv$$2$ (\textit{mod} $5$) \textit{or} \textit{if}  $p$$\equiv$$1$ (\textit{mod} $5$), $q$$\equiv$$3$ (\textit{mod} $5$), \textit{then the generalized Fibonacci-Lucas quaternions} $G^{p,q}_{n}$ \textit{are zero divisors in the quaternion ring} $\mathbb{H}_{\mathbb{Z}_{5}}\left( -1 ,-1 \right)$ \textit{if and only if} $n$$\equiv$$3$ (\textit{mod} $5$).\\
vi) \textit{if} $5$ \textit{does not divide} $p$ \textit{and} $q$$\equiv$$0$ (\textit{mod} $5$), \textit{then the generalized Fibonacci-Lucas quaternions} $G^{p,q}_{n}$ \textit{are zero divisors in the quaternion ring} $\mathbb{H}_{\mathbb{Z}_{5}}\left( -1 ,-1 \right)$ \textit{if and only if} $n$$\equiv$$2$ (\textit{mod} $5$).\\
\smallskip\\
\textbf{Proof.} i) If $p$$\equiv$$0$ (mod $5$) according to Proposition 3.1 we obtain that $n\left(G^{p,q}_{n}\right)=\overline{0}$ in $\mathbb{Z}_{5},$ $\left(\forall\right) n$$\in$$\mathbb{N}.$ This implies that all generalized Fibonacci-Lucas quaternions $G^{p,q}_{n}$ are zero divisors in the quaternion ring $\mathbb{H}_{\mathbb{Z}_{5}}\left( -1 ,-1 \right).$\\
ii) Applying Proposition 3.1, we have: $n\left(G^{p,q}_{n}\right)=\overline{3p^{2}f_{2n+1}+6pq\left(f_{2n+1}+f_{2n+3}\right)}$ in $\mathbb{Z}_{5}.$ Using the recurrence of Fibonacci sequence, it results that:
$$n\left(G^{p,q}_{n}\right)=\overline{p\cdot\left[\left(3p+2q\right)f_{2n+1}+qf_{2n+2}\right]}\ in\; \mathbb{Z}_{5}.$$
The generalized Fibonacci-Lucas quaternions $G^{p,q}_{n}$ are zero divisors in the quaternion ring $\mathbb{H}_{\mathbb{Z}_{5}}\left( -1 ,-1 \right)$ if and only if $n\left(G^{p,q}_{n}\right)=\overline{0}$ (in $\mathbb{Z}_{5}$).
Since $5$ does not divide $p,$ the $n\left(G^{p,q}_{n}\right)=\overline{0}$ (in $\mathbb{Z}_{5}$) if and only if $\left(3p+2q\right)f_{2n+1}+qf_{2n+2}$$\equiv$$0$ (mod $5$).
Since $p$$\equiv$$q$ (mod $5$), we obtain that the congruence
$\left(3p+2q\right)f_{2n+1}+qf_{2n+2}$$\equiv0$ (mod $5$) is equivalent with $qf_{2n+2}$$\equiv0$ (mod $5$) if and only in $f_{2n+2}$$\equiv0$ (mod $5$). According to Proposition 1.5, the last congruence is equivalent with $2n+2$$\equiv0$ (mod $5$) if and only if $n$$\equiv4$ (mod $5$).\\
iii) If $5$ does not divide $p$, similarly as in ii) the generalized Fibonacci-Lucas quaternions $G^{p,q}_{n}$ are zero divisors in the quaternion ring $\mathbb{H}_{\mathbb{Z}_{5}}\left( -1 ,-1 \right)$ if and only if $\left(3p+2q\right)f_{2n+1}+qf_{2n+2}$$\equiv0$ (mod $5$). Using the fact that $p+q$$\equiv$$0$ (mod $5$) and the recurrence of Fibonacci sequence, the last congruence is equivalent with $qf_{2n}$$\equiv0$ (mod $5$) if and only if
$f_{2n}$$\equiv0$ (mod $5$) if and only if $n$$\equiv 0$ (mod $5$) (according to Proposition 1.5).\\
iv) Similarly as in ii), the generalized Fibonacci-Lucas quaternions $G^{p,q}_{n}$ are zero divisors in the quaternion ring $\mathbb{H}_{\mathbb{Z}_{5}}\left( -1 ,-1 \right)$ if and only if $\left(3p+2q\right)f_{2n+1}+q f_{2n+2}$$\equiv0$ (mod $5$). We denote 
$E=\left(3p+2q\right)f_{2n+1}+q f_{2n+2}.$ We have:\\
$$E\equiv\left\{
\begin{array}{c}
2f_{2n+1}+2f_{2n+2}\; (mod\; 5),\ if\ p\equiv 1\; (mod\; 5),\ q\equiv 2\; (mod\; 5)\\
f_{2n+1}+f_{2n+2}\; (mod\; 5),\ if\ p\equiv 3\; (mod\; 5),\ q\equiv 1\; (mod\; 5)\\
4f_{2n+1}+4f_{2n+2}\; (mod\; 5),\ if\ p\equiv 2 \;(mod\; 5),\ q\equiv 4 \;(mod\; 5)\\
3f_{2n+1}+3f_{2n+2}\; (mod\; 5),\ if\ p\equiv 4 \;(mod\; 5),\ q\equiv 3\; (mod\; 5).\\

\end{array}%
\right.$$
This is equivalent with:\\
$$E\equiv\left\{
\begin{array}{c}
2f_{2n+3}\; (mod\; 5),\ if\ p\equiv 1\; (mod\; 5),\ q\equiv 2\; (mod\; 5)\\
f_{2n+3}\; (mod\; 5),\ if\ p\equiv 3\; (mod\; 5),\ q\equiv 1\; (mod\; 5)\\
4f_{2n+3}\; (mod\; 5),\ if\ p\equiv 2 \;(mod\; 5),\ q\equiv 4 \;(mod\; 5)\\
3f_{2n+3}\; (mod\; 5),\ if\ p\equiv 4 \;(mod\; 5),\ q\equiv 3\; (mod\; 5).\\

\end{array}%
\right.$$
For to have $E\equiv$ $0$ (mod $5$) it is necessary and sufficient to have $f_{2n+3}\equiv$$0$ (mod $5$). According to Proposition 1.5, the last congruence is equivalent with $2n+3$$\equiv0$ (mod $5$) if and only if $n$$\equiv1$ (mod $5$).\\
v) Similarly as in iv), we have:\\
$$E\equiv\left\{
\begin{array}{c}
3f_{2n+1}+f_{2n+2}\; (mod\; 5),\ if\ p\equiv 2\; (mod\; 5),\ q\equiv 1\; (mod\; 5)\\
2f_{2n+1}+4f_{2n+2}\; (mod\; 5),\ if\ p\equiv 3\; (mod\; 5),\ q\equiv 4\; (mod\; 5)\\
f_{2n+1}+2f_{2n+2}\; (mod\; 5),\ if\ p\equiv 4 \;(mod\; 5),\ q\equiv 2 \;(mod\; 5)\\
4f_{2n+1}+3f_{2n+2}\; (mod\; 5),\ if\ p\equiv 1 \;(mod\; 5),\ q\equiv 3\; (mod\; 5).\\

\end{array}%
\right.$$
Using Fibonacci recurrence, we obtain:
$$E\equiv\left\{
\begin{array}{c}
f_{2n}+4f_{2n+1}\; (mod\; 5),\ if\ p\equiv 2\; (mod\; 5),\ q\equiv 1\; (mod\; 5)\\
2f_{2n+3}+2f_{2n+2}\; (mod\; 5),\ if\ p\equiv 3\; (mod\; 5),\ q\equiv 4\; (mod\; 5)\\
f_{2n+3}+f_{2n+2}\; (mod\; 5),\ if\ p\equiv 4 \;(mod\; 5),\ q\equiv 2 \;(mod\; 5)\\
-f_{2n+1}-2f_{2n+2}\; (mod\; 5),\ if\ p\equiv 1 \;(mod\; 5),\ q\equiv 3\; (mod\; 5).\\

\end{array}%
\right.$$
This is equivalent with:\\
$$E\equiv\left\{
\begin{array}{c}
4f_{2n-1}\; (mod\; 5),\ if\ p\equiv 2\; (mod\; 5),\ q\equiv 1\; (mod\; 5)\\
2f_{2n+4}\; (mod\; 5),\ if\ p\equiv 3\; (mod\; 5),\ q\equiv 4\; (mod\; 5)\\
f_{2n+4}\; (mod\; 5),\ if\ p\equiv 4 \;(mod\; 5),\ q\equiv 2 \;(mod\; 5)\\
-f_{2n+4}\; (mod\; 5),\ if\ p\equiv 1 \;(mod\; 5),\ q\equiv 3\; (mod\; 5).\\

\end{array}%
\right.$$
Similarly as in ii), the generalized Fibonacci-Lucas quaternions $G^{p,q}_{n}$ are zero divisors in the quaternion ring $\mathbb{H}_{\mathbb{Z}_{5}}\left( -1 ,-1 \right)$ if and only if $E$$\equiv0$ (mod $5$). This is equivalent with $f_{2n-1}\equiv$$0$ (mod $5$), respectively $f_{2n+4}\equiv$$0$ (mod $5$). According to Proposition 1.5, each of these two congruences is equivalent to $n$$\equiv$$3$ (mod $5$).\\
vi) We get $E$$\equiv$$3pf_{2n+1}$ (mod $5$). For to have $E\equiv$$0$ (mod $5$) it is necessary and sufficient to have $f_{2n+1}\equiv$$0$ (mod $5$). According to Proposition 1.5, the last congruence is equivalent with $2n+1$$\equiv0$ (mod $5$) if and only if $n$$\equiv2$ (mod $5$).\\
\smallskip\\
Now, we consider $r=p\neq q,$ $p>5$ and we obtain the following result:\\
\smallskip\\
\textbf{Theorem 3.1.} \textit{Let} $p,q$ \textit{be two prime positive
integers}, $p\neq q,$ $p>5$ \textit{and let} $n$ \textit{be an arbitrary positive integer. Then, the generalized Fibonacci-Lucas quaternions} $G^{p,q}_{n}$ \textit{are zero divisors in the quaternion ring} $\mathbb{H}_{\mathbb{Z}_{p}}\left( -1 ,-1 \right)$ \textit{if and only if} $n=\frac{\left(2l+1\right)\cdot z\left(p\right)-1}{2}-1,$ \textit{where} $l$$\in$$\mathbb{N}.$\\
\smallskip\\
\textbf{Proof.} Applying Proposition 3.1 and the fact that $r=p,$ it results that
$\textbf{n}\left(G^{p,q}_{n}\right)=\overline{15q^{2}f_{2n+3}}$ in $\mathbb{Z}_{p}.$ Since $p,q$ are primes, $p\neq q,$ $p>5,$ we obtain that $\textbf{n}\left(G^{p,q}_{n}\right)=\overline{0}$ (in $\mathbb{Z}_{p}$) if and only if $p|f_{2n+3}.$ Similarly with the proof of Theorem 2.2 we obtain that $p|f_{2n+3}$ is equivalent with $n=\frac{\left(2l+1\right)\cdot z\left(p\right)-1}{2}-1,$ where $l$$\in$$\mathbb{N}.$\\
\smallskip\\
Now we consider $r=q\neq p,$ $q>5$ and we obtain the following result:\\
\smallskip\\
\textbf{Theorem 3.2.} \textit{Let} $p,q$ \textit{be two prime positive
integers}, $p\neq q,$ $q>5$ \textit{and let} $n$ \textit{be an arbitrary positive integer. Then, the}
\textit{generalized Fibonacci-Lucas quaternions} $G^{p,q}_{n}$ \textit{are zero divisors in the quaternion ring} $\mathbb{H}_{\mathbb{Z}_{q}}\left( -1 ,-1 \right)$ \textit{if and only if} $n=\frac{\left(2l+1\right)\cdot z\left(q\right)-1}{2},$ \textit{where} $l$$\in$$\mathbb{N}.$\\
\smallskip\\
\textbf{Proof.} Applying Proposition 3.1 and the fact that $r=q,$ it results that
$\textbf{n}\left(G^{p,q}_{n}\right)=\overline{3p^{2}f_{2n+1}}$ in $\mathbb{Z}_{q}.$ Since $p,q$ are primes, $p\neq q,$ $q>5,$ we obtain that $\textbf{n}\left(G^{p,q}_{n}\right)=\overline{0}$ (in $\mathbb{Z}_{q}$) if and only if $q|f_{2n+1}.$ Similarly with the proof of Theorem 2.2 we obtain that $q|f_{2n+1}$ is equivalent with $n=\frac{\left(2l+1\right)\cdot z\left(q\right)-1}{2},$ where $l$$\in$$\mathbb{N}.$\\
\begin{equation*}
\end{equation*}%
\textbf{References}
\begin{equation*}
\end{equation*}%
[Ak, Ko, To; 14] M. Akyigit, HH Kosal, M. Tosun, \textit{Fibonacci
Generalized Quaternions }, Adv. Appl. Clifford Algebras, vol. \textbf{24},
issue: 3 (2014), p. 631-641\newline
[Da, Dr; 70] D. E. Daykin, L. A. G. Dresel, \textit{Factorization of Fibonacci numbers}, The
Fibonacci Quarterly 8, No. 1 (1970).\newline
[Fl, Sa; 15] C. Flaut, D. Savin, \textit{Quaternion algebras and generalized Fibonacci-Lucas quaternions}, Adv. Appl. Clifford Algebras, vol. \textbf{25} issue 4 (2015), p. 853-862.\newline
[Gr, Mi, Ma; 15] J. M. Grau, C. Miguel and A. M. Oller-Marcen, \textit{On the Structure of Quaternion Rings Over} $\mathbb{Z}/n\mathbb{Z},$ Adv. Appl. Clifford Algebras, vol. \textbf{25}, issue: 4 (2015), p. 875-887\newline
[Ho; 61] A. F. Horadam, \textit{A Generalized Fibonacci Sequence}, Amer.
Math. Monthly, \textbf{68}(1961), p. 455-459.\newline
[Ho; 63] A. F. Horadam, \textit{Complex Fibonacci Numbers and Fibonacci
Quaternions}, Amer. Math. Monthly, \textbf{70}(1963), 289-291.\newline
[Lam; 04] T. Y. Lam, \textit{Introduction to Quadratic Forms over Fields,}
American Mathematical Society, 2004.\newline
[Mi, Se; 11] C. J. Miguel, R. Serodio, \textit{On the structure of quaternion rings over} $\mathbb{Z}_{p},$ International Journal of Algebra, vol.5, 2011, no.27, 1313-1325.\newline
[Po, Ki, Ke; 16]  E. Polatli, C. Kizilates, S. Kesim, \textit{On Split k-Fibonacci and k-Lucas Quaternions}, Adv. Appl. Clifford Algebras, vol. \textbf{26} issue 1 (2016), p. 353-362
\newline
[Su] Sun, Z. H. Congruences for Fibonacci numbers, http://www.hytc.cn/xsjl/szh\newline
[Su; 92] Sun, Z. H., Z.W. Sun. Fibonacci numbers and Fermat’s last theorem. Acta Arithemtica, Vol. 60, 1992, 371–388.\newline
[Ta; Ya; 15] D. Tasci, F. Yalcin, \textit{Fibonacci-}$p$ \textit{Quaternions}, Adv. Appl. Clifford Algebras, vol. \textbf{25} issue 1 (2015), p. 245-254.\newline
[Fib.] http://www.maths.surrey.ac.uk/hosted-sites/R.Knott/Fibonacci/fib.html%
\newline
[Y$\ddot{u}$, Ay; 16] S. Y$\ddot{u}$ce, F.T. Aydin, \textit{A New Aspect of Dual Fibonacci Quaternions,} accepted
in Adv. Appl. Clifford Algebras (http://link.springer.com/article/10.1007/\\
s00006-015-0619-9).\newline
\begin{equation*}
\end{equation*}

\bigskip

Diana SAVIN

{\small Faculty of Mathematics and Computer Science, }

{\small Ovidius University, }

{\small Bd. Mamaia 124, 900527, CONSTANTA, ROMANIA }

{\small http://www.univ-ovidius.ro/math/}

{\small e-mail: \ savin.diana@univ-ovidius.ro, \ dianet72@yahoo.com}\bigskip
\bigskip

\end{document}